\documentclass{article}
\usepackage{amsfonts}
\usepackage{amsmath}
\usepackage{enumerate}
\usepackage{amsthm}
\usepackage{hyperref}
\usepackage{cite}

\begin{document}

\newtheorem{Definition}{Definition}[section]
\newtheorem{Theorem}[Definition]{Theorem}
\newtheorem{Proposition}[Definition]{Proposition}
\newtheorem{Remark}[Definition]{Remark}
\newtheorem{Lemma}[Definition]{Lemma}
\newtheorem{Corollary}[Definition]{Corollary}

\numberwithin{equation}{section}

\title{Prescribed Webster scalar curvatures on compact pseudo-Hermitian manifolds}
\author{Yuxin Dong\footnote{Supported by NSFC grant No. 11771087, and LMNS, Fudan},\ Yibin Ren\footnote{Supported by NSFC grant No. 11801517},\ Weike Yu\footnote{Corresponding author}}
\date{}
\maketitle
\begin{abstract}
In this paper, we investigate the problem of prescribing Webster scalar curvatures on compact pseudo-Hermitian manifolds. In terms of the method of upper and lower solutions and the perturbation theory of self-adjoint operators, we can describe some sets of Webster scalar curvature functions which can be realized through pointwise CR conformal deformations and CR conformally equivalent deformations respectively from a given pseudo-Hermitian structure.
\end{abstract}

\section{Introduction}
In Riemannian geometry, the problem of finding a conformal metric on a compact Riemannian manifold with a prescribed scalar curvature has been investigated extensively (cf. \cite{[KW1], [KW2], [KW3], [Ou1], [Ou2], [Ta], [Ra], [Ho1], [CX]} and the references therein). Its special case that the candidate scalar curvature function is constant is the well-known Yamabe problem, which was settled down by a series of works due to Yamabe, Trudinger, Aubin, and Schoen (cf. \cite{[Yam], [Tr], [Au], [Sc]}).

The following problem is a CR analogue of prescribing scalar curvature problem: given any smooth function $\hat{\rho}$ on a compact strictly pseudoconvex CR manifold $M$ of real dimension $2n+1$ with contact form $\theta$, does there exist a contact form $\hat{\theta}$ CR conformal to $\theta$, that is, $\hat{\theta}=u^{\frac{2}{n}}\theta$ for some positive function $u$, such that its Webster scalar curvature $\text{Scal}_{\hat{\theta}}=\hat{\rho}$? It is equivalent to solving the following partial differential equation
\begin{align}
-(2+\frac{2}{n})\Delta_\theta u+\text{Scal}_\theta u=\hat{\rho}u^{1+\frac{2}{n}}\ \ \ \text{on}\ M\label{1.1.}
\end{align}
for $u>0$, where $\text{Scal}_\theta$ is the Webster scalar curvature of $(M,\theta)$. When $\hat{\rho}$ is constant, the above problem is referred to as CR Yamabe problem, which was solved by Jerison and Lee (cf. \cite{[JL1], [JL2]}), Gamara and Yacoub (cf. \cite{[Ga1], [GY]}). Another interesting special case for the prescribed Webster scalar curvature problem is to consider the domain manifold to be a CR sphere $S^{2n+1}$. Similar to the Riemannian case, this problem is not always solvable. Indeed, Cheng \cite{[Ch]} gave a Kazdan-Warner type necessary condition for the solution $u$ and the prescribed function $\hat{\rho}$. Besides, in \cite{[FU], [MU], [HK], [SG], [RG], [CPY], [Ho3], [Ho4], [Ho5]}, if $\hat{\rho}$ satisfies suitable conditions, some existence results were established for the prescribed Webster scalar curvature problem on $S^{2n+1}$ by means of variational, topological, perturbation methods, Webster scalar curvature flow, or the theory of critical points, etc. In \cite{[Ga2], [CEG], [CAY1], [CAY2], [Yac], [GAG]}, the authors investigated the problem on strictly pseudoconvex spherical CR manifolds. In \cite{[Ho2]}, using geometric flow, Ho proved that any negative smooth function $\hat{\rho}$ can be prescribed as the Webster scalar curvature in the CR conformal class, provided that $\dim M=3$ and the CR Yamabe invariant of M is negative. In \cite{[NZ]}, the authors studied the prescribed Webster scalar curvature problem on a pseudo-Hermitian manifold in arbitrary CR dimension with negative CR Yamabe invariant. Using variational techniques, they established several non-existence, existence, and multiplicity results when the function $\hat{\rho}$ is sign-changing.

In this paper, we investigate the prescribed Webster scalar curvature problem on a compact strictly pseudoconvex CR manifold $M$, following the original approaches in \cite{[KW1],[KW2]}, but adjusting their related arguments to subelliptic version. In this way, we are able to establish several existence results for the problem on a compact strictly pseudoconvex CR manifold in arbitrary CR dimension. To state our main results, let us introduce some notations. Given a compact strictly pseudoconvex CR manifold $(M^{2n+1},H,J,\theta)$ (also called the pseudo-Hermitian manifold, see Section \ref{section2} ), where $(H,J)$ is a CR structure of type $(n,1)$ and $\theta$ is a pseudo-Hermitian structure with positive Levi form, let $PC(\theta)$ denote the set of smooth functions on $M$ that are the Webster scalar curvatures of pseudo-Hermitian structures $\hat{\theta}$ in the CR conformal class $[\theta]=\{u\theta: 0<u\in C^\infty(M)\}$. In other words, $PC(\theta)$ is the set of smooth functions $\hat{\rho}$ for which one can find a positive solution of \eqref{1.1.}. Let $Y_M(\theta)$ be the CR Yamabe constant (see \eqref{2.14} or \eqref{2.15}) and $\lambda_1$ be the first eigenvalue of operator $L=-(2+\frac{2}{n})\Delta_\theta+\text{Scal}_\theta$ (also see \eqref{2.12}). Using the method of upper and lower solutions on CR manifolds, we obtain the following conclusions.
\begin{Theorem}\label{theorem3.10}
Let $(M^{2n+1},H,J,\theta)$ be a compact pseudo-Hermitian manifold. Then 
\begin{enumerate}[(A)]
\item The following statements are equivalent:\\
$(a1)$ $\lambda_1<0$.\\
$(a2)$ $\{\hat{\rho}\in C^\infty(M):\hat{\rho}<0\} \subset PC(\theta).$\\
$(a3)$ $\{\hat{\rho}\in C^\infty(M):\hat{\rho}<0\} \cap PC(\theta)\neq \emptyset.$\\
$(a4)$ $Y_M(\theta)<0.$
\item The following statements are equivalent:\\
$(b1)$ $\lambda_1=0$\\
$(b2)$ $0\in PC(\theta).$\\
$(b3)$ $ Y_M(\theta)=0.$
\item The following statements are equivalent:\\
$(c1)$ $\lambda_1>0.$\\
$(c2)$ $\{\hat{\rho}\in C^\infty(M):\hat{\rho}>0\} \cap PC(\theta)\neq \emptyset.$\\
$(c3)$ $Y_M(\theta)>0.$
\end{enumerate}
\end{Theorem}
Note that the results in (B) of Theorem 1.1  essentially belong to \cite{[JL1]} as a special case of CR Yamabe problem, and we state these results here for relative completeness, while the prescribed Webster scalar curvature problem for more general nonvanishing function $\hat{\rho}$ in the case of $\lambda_1=0$ is still an unsolved problem.

Since in general not all smooth functions $\hat{\rho}$ can be realized as the Webster scalar curvature of some pseudo-Hermitian structure $\hat{\theta}$ pointwise CR conformal to $\theta$, i.e., $\hat{\theta}\in [\theta]$ (cf. \cite{[Ch]}, \cite{[NZ]}), we will try to enhance the possibility of realizing $\hat{\rho}$ as the prescribed Webster scalar curvature by relaxing the desired pseudo-Hermitian structure $(\hat{H},\hat{J}, \hat{\theta})$ to be CR conformally equivalent to $(H,J,\theta)$, i.e., there is a map $\Phi\in \text{Diff}(M)$ such that $\Phi^*\hat{\theta}\in [\theta],\hat{H}=d\Phi(H),\hat{J}=d\Phi\circ J\circ (d\Phi)^{-1}$. For simplicity, let $CE(\theta)$ denote the set of smooth functions on $M$ which are the Webster scalar curvatures of $(\hat{H},\hat{J}, \hat{\theta})$ CR conformally equivalent to $(H,J,\theta)$. In other words, $CE(\theta)$ is the set of smooth functions $\hat{\rho}$ for which one can find a map $\Phi\in\text{Diff}(M)$ such that
\begin{align}
-(2+\frac{2}{n})\Delta_\theta u+\text{Scal}_\theta u=(\hat{\rho}\circ \Phi) u^{1+\frac{2}{n}}\ \ \ \text{on}\ M
\end{align}
admits a positive solution. By the inverse function theorem and perturbation methods in our cases, we obtain
\begin{Theorem}\label{theorem1.2}
Let $(M^{2n+1},H,J,\theta)$ be a compact pseudo-Hermitian manifold.
\begin{enumerate}[(1)]
\item If $\lambda_1<0$, then $CE(\theta)=\{\hat{\rho}\in C^\infty(M):\ \hat{\rho}<0\ \text{somewhere}\}$.
\item If $\lambda_1=0$, then $CE(\theta)=\{\hat{\rho}\in C^\infty(M):\ \hat{\rho}\ \text{changes\ sign\ on\ M}\}\cup \{0\}$.
\item If $\lambda_1>0$, then $CE(\theta)=\{\hat{\rho}\in C^\infty(M):\ \hat{\rho}>0\ \text{somewhere}\}$.
\end{enumerate}
\end{Theorem}
In particular, if $\hat{\rho}$ is a smooth function on $M$ with changing sign, then it belongs to $CE(\theta)$ regardless of the sign of $\lambda_1$. In other words, any smooth function with changing sign can be realized as some Webster scalar curvature.
 \begin{Corollary}
Let $(M^{2n+1},H,J,\theta)$ be a compact pseudo-Hermitian manifold. If $\hat{\rho}\in C^\infty(M)$ and it changes sign on $M$, then there exists a structure $(\hat{H},\hat{J},\hat{\theta})$ on $M$ such that $(M,\hat{H},\hat{J},\hat{\theta})$ is a pseudo-Hermitian manifold with the Webster scalar curvature $\hat{\rho}$ and is CR conformally equivalent to $(M,H,J,\theta)$.
 \end{Corollary}

\section{Preliminaries}\label{section2}
In this section, we will introduce the notions and notations of pseudo-Hermitian geometry (cf. \cite{[DT]}).

Let $M$ be an orientable real smooth manifold with $\dim_{\mathbb{R}} M=2n+1$. A CR structure on $M$ is a complex subbundle $T_{1,0}M$ of complex rank $n$ of the complexified tangent bundle $TM\otimes\mathbb{C}$ satisfying
\begin{align}
T_{1,0}M\cap T_{0,1}M=\{0\}, \ \ [\Gamma(T_{1,0}M),\Gamma(T_{1,0}M)]\subseteq\Gamma(T_{1,0}M)\label{2.1}
\end{align}
where $T_{0,1}M=\overline{T_{1,0}M}$. The complex subbundle $T_{1,0}M$ corresponds to a real rank $2n$ subbundle of $TM$:
\begin{align}
H=Re\{T_{1,0}M\oplus T_{0,1}M\},
\end{align}
which is called the Levi distribution. Clearly, it carries a natural complex structure $J$ defined as
\begin{align}
J(X+\bar{X})=\sqrt{-1}(X-\bar{X})
\end{align}
for any $X\in T_{1,0}M$. Equivalently, the CR structure may be described by the pair $(H,J)$. Let $(M,H,J)$ and $(\tilde{M},\tilde{H},\tilde{J})$ be two CR manifolds. A smooth map $f: (M,H,J) \rightarrow (\tilde{M},\tilde{H},\tilde{J})$ is called a CR map if it satisfies
\begin{align}
df(H)\subset \tilde{H},\ \ \ \ df\circ J=\tilde{J}\circ df\ \text{on}\ H.
\end{align}
Furthermore, $f$ is said to be a CR isomorphism if it is a $C^\infty$ diffeomorphism and a CR map.

 Since both $M$ and $H$ are orientable, there is a global nowhere vanishing 1-form $\theta$ with $H=\ker \theta$, which is called a pseudo-Hermitian structure on $M$. The corresponding Levi form is defined as
\begin{align}
L_\theta(X,Y)=d\theta(X,JY)
\end{align}
for any $X,Y\in H$. The integrability assumption of $T_{1,0}M$ implies $L_{\theta}$ is $J$-invariant and symmetric. If the CR manifold $M$ admits a pseudo-Hermitian structure $\theta$ such that $L_\theta$ is positive definite, then $(M,H,J)$ is said to be strictly pseudoconvex. Henceforth we will assume that $(M,H,J)$ is a strictly pseudoconvex CR manifold and $\theta$ is a pseudo-Hermitian structure with positive Levi form. The quadruple $(M^{2n+1}, H, J, \theta)$ is referred to as a pseudo-Hermitian manifold. 

Let $\theta, \hat{\theta}$ be two pseudo-Hermitian structures on the CR manifold $(M,H,J)$, whose Levi forms are positive definite. Since $\dim_{\mathbb{R}}TM/H=1$, $\theta, \hat{\theta}$ are related by
 \begin{align}
 \hat{\theta}=u\theta \label{2.5.}
 \end{align}
 for some nowhere vanishing function $u\in C^\infty(M)$. Applying the exterior differentiation operator $d$ to \eqref{2.5.}, we get
 \begin{align}
 L_{\hat{\theta}}=uL_\theta.
 \end{align}
 Since both $L_\theta$ and $L_{\hat{\theta}}$ are positive definite, we see that $u$ is positive everywhere. Given a CR structure $(H,J)$, then the set of all its pseudo-Hermitian structures with positive Levi form is exactly
 \begin{align}
 [\theta]=\{u\theta:\ 0<u\in C^\infty(M)\},
 \end{align}
 where $\theta$ is one pseudo-Hermitian structure of $(H,J)$ with positive Levi form.
 A property on a CR manifold $(M,\theta)$ is said to be CR invariant if it is invariant for all pseudo-Hermitian structures in $[\theta]$. 

On a pseudo-Hermitian manifold $(M^{2n+1}, H, J, \theta)$, there is a unique globally defined nowhere vanishing tangent vector field $\xi$ on M such that \begin{align}
\theta(\xi)=1,\ \ \ \ d\theta(\xi,\cdot)=0,
\end{align}
which is usually called the Reeb vector field. Hence, we have a splitting of the tangent bundle 
\begin{align}
TM=H\oplus \mathbb{R}\xi,
\end{align}
which leads to a natural projection $\pi_H:TM\rightarrow H$ and a Riemannian metric on $M$ (the Webster metric)
\begin{align}
g_\theta=\pi^*_HL_\theta+\theta\otimes\theta,
\end{align}
where $(\pi^*_HL_\theta)(X,Y)=L_{\theta}(\pi_HX, \pi_HY)$ for $X,Y\in TM$. On a pseudo-Hermitian manifold, there is a unique linear connection $\nabla$ called Tanaka-Webster connection preserving the CR structure and Webster metric (cf. Theorem 1.3 of \cite{[DT]} ). For a smooth function $u$ on $M$, one can define the sub-Laplacian of $u$ as the divergence of horizontal gradient:
\begin{align}
\Delta_\theta u=\text{div} (\nabla^H u),
\end{align}
where $\nabla^H u=\pi_H \nabla u$. Then the integration by parts yields
\begin{align}
\int_M (\Delta_\theta u)v \Psi^{\theta}=-\int_M L_\theta(\nabla^H u, \nabla^H v)\Psi^{\theta},
\end{align}
for $u,v \in C^2(M)$ with compact support, where $\Psi^{\theta}=\theta\wedge (d\theta)^n$ is a volume form of $(M^{2n+1},H,J,\theta)$.

The curvature theory of Tanaka-Webster connection was developed in \cite{[We]} (cf. also \cite{[DT]}). In particular, Webster defined a scalar curvature associated with a pseudo-Hermitian structure $\theta$, which is referred to as Webster scalar curvature in literature.
For a pseudo-Hermitian manifold $(M^{2n+1},H,J,\theta)$, we say that a pseudo-Hermitian structure $\hat{\theta}$ on $M$ is pointwise CR conformal to $\theta$ if $\hat{\theta}=u^{\frac{2}{n}}\theta$ for some positive function $u\in C^\infty(M)$. The pseudo-Hermitian manifold $(M^{2n+1},H,J,\hat{\theta})$ is said to be a pointwise CR conformal deformation of $(M^{2n+1},H,J,\theta)$. Furthermore, from \cite{[Le]} and \cite{[JL1]}, the Webster scalar curvatures of $\theta$ and $\hat{\theta}$ have the following relationship:
\begin{align}
-b_n\Delta_\theta u+\rho u=\hat{\rho}u^a\label{2.10.}
\end{align}
where $\hat{\rho}$ is the Webster scalar curvature of $(M,H,J,\hat{\theta})$, and $a=1+\frac{2}{n}$, $b_n=2+\frac{2}{n}$. For convenience, let $PC(\theta)$ denote the set of $C^\infty(M)$ functions which are Webster scalar curvatures of all $\hat{\theta}\in [\theta]$. In other words, $PC(\theta)$ is the set of $C^\infty(M)$ functions for which one can find a positive solution of \eqref{2.10.}.

Now we assume that $(M^{2n+1},H,J,\theta)$ is a compact pseudo-Hermitian manifold. Set
\begin{align}
L=-b_n\Delta_\theta+\rho\label{2.15...}
\end{align}
and let $\lambda_1$ be the first eigenvalue of the operator $L$, that is 
\begin{align}
\lambda_1=\inf_{u\in S^2_1(M)-\{0\}}\frac{\int_M (b_n|\nabla^Hu|_\theta^2+\rho u^2)\Psi^\theta}{\int_M u^2\Psi^\theta},
\label{2.12}
\end{align}
where $ S^2_1(M)$ is the Folland-Stein space (cf. \cite{[FS]}, \cite{[RS]}), the norm $|\cdot|_\theta$ is induced by $g_\theta$. If $\psi$ is an eigenfunction corresponding to $\lambda_1$, then $L\psi=\lambda_1\psi$. Note that $\psi$ is $C^\infty$ and nowhere vanishing (cf. \cite{[Wa]}), so we may assume that $\psi>0$ and thus
\begin{align}
\lambda_1=\inf_{0<u\in C^\infty(M)}\frac{\int_M (b_n|\nabla^Hu|_\theta^2+\rho u^2)\Psi^\theta}{\int_M u^2\Psi^\theta}.
\end{align}
Recall that CR Yamabe constant is given by
\begin{align}
Y_M(\theta)&=\inf_{0<u\in C^\infty(M)}\frac{\int_M (b_n|\nabla^Hu|_\theta^2+\rho u^2)\Psi^\theta}{(\int_M u^{2+\frac{2}{n}}\Psi^\theta)^{\frac{n}{n+1}}}\label{2.14}\\
        &=\inf_{\hat{\theta}\in [\theta]} \frac{\int_M \hat{\rho} \Psi^{\hat{\theta}} }{(\int_M \Psi^{\hat{\theta}} )^{\frac{n}{n+1}}}\label{2.15}
\end{align}
which is a CR invariant.

Given a pseudo-Hermitian manifold $(M^{2n+1},H,J,\theta)$, we say that the structures $(\hat{H},\hat{J},\hat{\theta})$ is CR conformally equivalent to $(H,J,\theta)$ if there is a map $\Phi\in \text{Diff}(M)$ and $0<u\in C^\infty(M)$ such that
\begin{align}
\Phi^*\hat{\theta}=u^{\frac{2}{n}}\theta,\ \ \ \hat{H}=d\Phi(H),\ \ \ \hat{J}=d\Phi\circ J\circ (d\Phi)^{-1}.\label{2.19..}
\end{align}
Clearly, $\hat{J}$ is a complex structure on $\hat{H}$ and $\Phi: (M,H,J,\theta)\rightarrow (M, \hat{H}, \hat{J}, \hat{\theta})$ is a CR isomorphism, where the pseudo-Hermitian manifold $(M, \hat{H}, \hat{J}, \hat{\theta})$ is called a CR conformally equivalent deformation of $(M,H,J,\theta)$. Furthermore, Webster scalar curvatures have the following relationship:
\begin{align}
-b_n\Delta_\theta u+\rho u=(\hat{\rho}\circ \Phi) u^a,\label{2.21.}
\end{align}
where $\hat{\rho}$ is the Webster scalar curvature of $(M, \hat{H}, \hat{J}, \hat{\theta})$, and $a=1+\frac{2}{n}$, $b_n=2+\frac{2}{n}$. Similarly,  let $CE(\theta)$ denote the set of $C^\infty(M)$ functions which are the Webster scalar curvatures of $(M, \hat{H}, \hat{J}, \hat{\theta})$. In other words, $CE(\theta)$ is the set of $C^\infty(M)$ functions $\hat{\rho}$ for which one can find a map $\Phi\in\text{Diff}(M)$ such that \eqref{2.21.} admits a positive solution on $M$. Clearly, $PC(\theta)$ is a subset of $CE(\theta)$.

At the end of this section, we recall the Folland-Stein spaces on the pseudo-Hermitian manifold $(M^{2n+1},H,J,\theta)$ briefly (cf. \cite{[FS]}, \cite{[RS]}), which are the generalized Sobolev spaces compatible to the CR structure $(H,J)$. Let $\{X_\alpha\}_{\alpha=1}^{2n}$ be a local $G_\theta$-orthonormal real frame of $H$ defined on an open subset $U\subset M$. For any $k\in\mathbb{N}_+$ and $1<p<+\infty$, the Folland-Stein spaces on $U$ is defined by
\begin{align}
S^p_k(U)=\{f\in L^p(U) : X_{i_1}X_{i_2}\dots X_{i_s}f\in L^p(U), s\leq k, X_{i_j}\in \{X_{\alpha}\}\}
\end{align}
with the norms
\begin{align}
\|f\|_{S^p_k(U)}=\|f\|_{L^p(U)}+\sum_{1\leq s\leq k} \|X_{i_1}X_{i_2}\dots X_{i_s}f\|_{L^p(U)}
\end{align}
where the $L^p$-norm of $f$ is defined by $\|f\|_{L^p(U)}=\left(\int_U |f|^p \Psi^\theta\right)^{\frac{1}{p}}$. By the partition of unity, we can also define $S^p_k(\Omega)$ and $S^p_k(M)$, where $\Omega$ is any open subset of $M$.

\section{ Pointwise CR conformal deformations with prescribed Webster scalar curvature}
\label{section3}
In this section, we will investigate the set $PC(\theta)$ on a compact pseudo-Hermitian manifold $(M^{2n+1},H,J,\theta)$ when $\lambda_1<0$, $\lambda_1=0$ and $\lambda_1>0$ respectively. Firstly, we consider the case $\lambda_1<0$. For this case, we will use the method of upper and lower solutions on pseudo-Hermitian manifolds. For this purpose, we need the following existence and comparison results.

\begin{Lemma}\label{lemma4.1}
Let $(M^{2n+1},H,J,\theta)$ be a compact pseudo-Hermitian manifold. Let $L_1=-\Delta_\theta+f$, where $f\in C^\infty(M)$. Then
\begin{enumerate}[(1)]
\item If $f>0$, then $L_1: S^2_2(M)\rightarrow L^2(M)$ is invertible.
\item  The equation $L_1v=g$ has a weak solution if and only if $\langle g,w\rangle_{L^2}=0$ for any solution $w$ of $L_1w=0$.
\item If  $f>0$ and $v$ is a $S^2_1(M)$ function with
\begin{align}
L_1v\geq 0,
\end{align}
then $v\geq 0$.

\end{enumerate}
\end{Lemma}
\proof
(1) Let us show that for any $g\in L^2(M)$, there is a unique solution $u\in S^2_2(M)$ such that $L_1u=g$. Here we can relax the requirement $u\in S^2_2(M)$ to $u\in S^2_1(M)$, because of the $L^2$ interior regularity result for $\Delta_\theta$. Set
\begin{align}
(u,v)=\langle L_1u,v\rangle_{L^2}=\int_M \left(\nabla^Hu\cdot\nabla^H v+ fuv\right) \Psi^\theta,
\end{align}
where $\cdot$ is the inner product induced by the Webster metric $g_\theta$ and $u,v\in S^2_1(M)$.
By a simple computation, we obtain 
\begin{align}
(u,u)\geq c_1\|u\|_{S^2_1(M)}^2\label{3.3...}
\end{align}
and
\begin{align}
(u,v)\leq c_2\|u\|_{S^2_1(M)}\|v\|_{S^2_1(M)}
\end{align}
where $c_1,c_2$ are two positive constants. Therefore, the space $S^2_1(M)$ with the inner product $(\cdot,\cdot)$ is a Hilbert space. By Cauchy-Schwarz inequality and \eqref{3.3...}, we have
\begin{align}
|\langle g,v\rangle_{L^2}|\leq \|g\|_{L^2}\|v\|_{L^2}\leq c_1^{-\frac{1}{2}}\|g\|_{L^2}(v,v)^{\frac{1}{2}}
\end{align}
for any $v\in S^2_1(M)$, which implies $\langle g,v\rangle_{L^2}$ is a bounded linear functional of $v\in S^2_1(M)$. Applying the Riesz representation theorem, there is a unique $u\in S^2_1(M)$ such that
\begin{align}
\langle g,v\rangle_{L^2}=(u,v)=\langle L_1u,v\rangle_{L^2}
\end{align}
for any $v\in  S^2_1(M)$.

(2) Since $M$ is compact and $f\in C^\infty(M)$, there is a positive constant $\lambda>0$ such that $f+\lambda>0$ on $M$. In view of the part (1) of this lemma, the inverse operator $L_2=(L_1+\lambda)^{-1}:\ L^2(M)\rightarrow S^2_2(M)$ exists. Using $S^2_1(M)\subset \subset L^2(M)$ (cf. Theorem 3.15 of \cite{[DT]}, \cite{[Da]}) yields that $L_2: L^2(M)\rightarrow L^2(M)$ is a completely continuous. The equation $L_1v=g$ is equivalent to $v-\lambda L_2v=L_2g$. Applying the Fredholm-Riesz-Schauder theory (cf. \cite{[BJS]}) and the facts $L_1^*=L_1, L_2^*=L_2$, we get that $v-\lambda L_2v=L_2g$ has a weak solution if and only if $\langle L_2g,w\rangle_{L^2}=0$ where $w$ satisfies $w-\lambda L_2w=0$ which is equivalent to $L_1w=0$. From
\begin{align}
\langle g,w\rangle_{L^2}=\langle L_2^{-1}L_2g,w\rangle_{L^2}=\langle L_2g,L_2^{-1}w\rangle_{L^2}=\lambda \langle L_2g, w\rangle_{L^2},
\end{align}
it follows that $L_1v=g$ has a weak solution if and only if $\langle g, w\rangle_{L^2}=0$ for any solution $w$ of $L_1w=0$.

(3) Since $v\in S^2_1(M)$, $v_-=\min\{v,0\}\in S^2_1(M)$. Taking $-v_-$ as a test function of $L_1v\geq0$ yields
\begin{align}
\int_M |\nabla^H v_-|_\theta^2 \Psi^\theta\leq -\int_M f(v_-)^2\Psi^\theta,\label{3.7...}
\end{align}
which implies $v_{-}=0$ since $f>0$. Hence, $v\geq 0$.
\qed

Using the above lemma, we obtain the following result, which is a pseudo-Hermitian version of Lemma 2.6 of \cite{[KW1]}.

\begin{Lemma}\label{lemma4.2}
Let $(M^{2n+1},H,J,\theta)$ be a compact pseudo-Hermitian manifold. Assume that $f(x,u)\in C^\infty(M\times \mathbb{R})$. If there are two functions $u_+,u_-\in C^0(M)\cap S^2_1(M)$ satisfying
\begin{equation}
-\Delta_\theta u_++f(x, u_+)\geq 0\ \ \ \ \text{in}\ M,\label{4.1}
\end{equation}
\begin{equation}
-\Delta_\theta u_-+f(x, u_-)\leq 0\ \ \ \ \text{in}\ M\label{4.2}\\
\end{equation}
\begin{equation}
u_+\geq u_- \ \ \ \ \text{in}\ M,\label{3.11...}
\end{equation}
then there exists a function $u\in C^\infty(M)$ such that
\begin{equation}
-\Delta_\theta u+f(x, u)=0\ \ \ \ \text{in}\ M,\label{4.4}
\end{equation}
\begin{equation}
u_-\leq u\leq u_+ \ \ \ \  \text{in}\ M.\label{4.5}
\end{equation}
Here $u_+$ and $u_-$ are called the upper and lower solutions of \eqref{4.4} respectively.
\end{Lemma}
\proof
Set $A_1=\min_M u_-$, $A_2=\max_M u_+$ and $I=[A_1, A_2]$. Since $f(x,u)\in C^\infty(M\times \mathbb{R})$, there exists a constant $\lambda>0$ such that $\tilde{f}(x,u)=-f(x,u)+\lambda u$ is increasing with respect to $u\in I$ for any fixed $x\in M$. In order to find a solution of \eqref{4.4} and \eqref{4.5}, we consider the sequence $\{u_k\}$ defined for $k\geq 1$ by
\begin{equation}
\left\{ 
\begin{aligned}
&(-\Delta_\theta+\lambda) u_k=\tilde{f}(x, u_{k-1}) \\
&u_0=u_-.
\end{aligned}
\right.
\label{4.6}
\end{equation}
By \eqref{4.1}, \eqref{4.2}, \eqref{3.11...} and \eqref{4.6}, we have
\begin{equation}
(-\Delta_\theta+\lambda)(u_+-u_1)\geq 0,
\end{equation}
\begin{equation}
(-\Delta_\theta+\lambda) (u_1-u_-)\geq 0.\label{4.8}
\end{equation}
According to Lemma \ref{lemma4.1}, 
\begin{align}
u_-\leq u_1\leq u_+.
\end{align}
Iterating the above procedure yields
\begin{align}
u_-\leq u_1\leq u_2\leq \ldots\leq u_+.
\end{align}
Set $u=\lim_{k\rightarrow\infty}u_k$, then $u$ is a weak solution of \eqref{4.4} and \eqref{4.5}. Since $\text{Im}\ u\subset I$ and $f(\cdot,\cdot)\in C^\infty(M\times \mathbb{R})$, we conclude that $f(x,u)\in L^p(M)$ with $p>2n+1$. By the regularity result for $\Delta_\theta$ (cf. Theorem 18 of \cite{[RS]}), we have $u\in S^p_2(M)$, and thus $f(x,u)\in S_2^p(M)$. Repeating the above argument, we obtain that $u\in S^p_{2k}(M)$ for any $k\in\mathbb{N}_+$. Therefore, $u\in C^\infty(M)$ because of $S^p_{2k}(M)\subset W^{k,p}(M) \subset C^{k-1}(M)$ for any $k\in\mathbb{N}_+$ (cf. Theorem 19.1 of \cite{[FS]}), where $W^{k,p}(M)$ is the classical Sobolev space.
\qed
\begin{Remark}
The authors of \cite{[NZ]} proved that when $f(x,u)=-b_n^{-1}(\rho u-\hat{\rho}u^a)$, the equation \eqref{4.4} admits a weak solution satisfying \eqref{4.5} if it has a weak lower solution $u_-$ and a weak upper solution $u_+$. 
\end{Remark}

In terms of Lemma \ref{lemma4.2}, and by a similar argument of \cite{[KW1]},  we obtain

\begin{Theorem}
\label{theorem3.4.}
Let $(M^{2n+1},H,J,\theta)$ be a compact pseudo-Hermitian manifold and $\hat{\rho}$ be a smooth negative function on $M$. Then $\hat{\rho}\in PC(\theta)$ if and only if $\lambda_1<0$.
\end{Theorem}
\proof
If $\hat{\rho}\in PC(\theta)$, then there is a positive function $u\in C^\infty(M)$ such that $u$ satisfies the prescribed Webster scalar curvature equation \eqref{2.10.}, that is, $Lu=\hat{\rho}u^a$, where $L$ is given by \eqref{2.15...} and $a=1+\frac{2}{n}$. Let $\psi$ be the positive eigenfunction associated with $\lambda_1$ of $L$. Then
\begin{align}
\lambda_1\langle \psi, u\rangle_{L^2}=\langle L\psi, u\rangle_{L^2}=\langle \psi, Lu\rangle_{L^2}=\langle \psi, \hat{\rho} u^a\rangle_{L^2}<0
\end{align}
from which it follows that $\lambda_1<0$. Conversely, If $\lambda_1<0$, then there is a positive solution $u\in C^\infty(M)$ such that $Lu=\hat{\rho}u^a$ by the existence of upper and lower solutions, hence $\hat{\rho}\in PC(\theta)$. Indeed, Let $u_+\equiv \alpha>0$ where $\alpha$ is a constant large enough so that
\begin{align}
Lu_+-\hat{\rho}u_+^a=\alpha(\rho-\hat{\rho}\alpha^{a-1})\geq 0.
\end{align}
On the other hand, let $u_-=\beta\psi$ where $\beta>0$ is so small that $u_-\leq \alpha\equiv u_+$ and $u_-\leq\left( \frac{\lambda_1}{\inf_M \hat{\rho}}\right)^{\frac{1}{a-1}}$. Then
\begin{align}
Lu_-=\beta L\psi=\lambda_1\beta\psi=\lambda_1u_-\leq  \hat{\rho} u_-^a.\label{3.16.}
\end{align}
Therefore, by Lemma \ref{lemma4.2}, there exists a smooth solution $u$ satisfying $Lu=\hat{\rho} u^a$ and $0<u_-\leq u\leq u_+$, i.e., $\hat{\rho}\in PC(\theta)$.
\qed
\begin{Remark}\label{remark3.5}
From the above proof, if $\lambda_1<0$, then there always exists a lower solution $u_-\in C^\infty(M)$ of \eqref{2.10.} with $0<u_-<u$, where $u$ is a given $C^\infty(M)$ function.
\end{Remark} 
\begin{Remark}
By a flow method, Ho \cite{[Ho2]} proved that every negative function $\hat{\rho}\in PC(\theta)$ if the CR Yamabe constant $Y_M(\theta)<0$ and $\dim_{\mathbb{R}} M=3$.
\end{Remark}
\begin{Remark}
In \cite{[NZ]}, authors gave the following results. 
\begin{enumerate}[(1)]
\item When $\hat{\rho}$ is a smooth nonpositive function on $(M^{2n+1},H,J,\theta)$ with $Y_M(\theta)<0$ such that the set $\{x\in M: \hat{\rho}(x)=0\}$ has positive measure, authors gave a necessary and sufficient condition for $\hat{\rho}\in PC(\theta)$.
\item If $\hat{\rho}$ is a smooth nonpositive function on $(M^{2n+1},H,J,\theta)$, then $\hat{\rho}$ is the Webster scalar curvature of at most one pseudo-Hermitian structure $\hat{\theta}\in [\theta]$.
\end{enumerate}
\end{Remark}

Making use of Lemma \ref{lemma4.2} again, we can establish the following property of the set $PC(\theta)$.

\begin{Proposition}\label{theorem3.7.}
Let $(M^{2n+1},H,J,\theta)$ be a compact pseudo-Hermitian manifold with $\lambda_1<0$. If $\hat{\rho}\in PC(\theta)$ and $\hat{\rho}_1\leq \hat{\rho}$, then $\hat{\rho}_1\in PC(\theta)$.
\end{Proposition}
\proof
To prove  $\hat{\rho}_1\in PC(\theta)$, we just need to find a positive solution of 
\begin{align}
-b_n\Delta_\theta u+\rho u=\hat{\rho}_1 u^a,\label{3.21.}
\end{align}
where $\rho$ is the Webster scalar curvature of $(M^{2n+1},H,J,\theta)$. We will use the method of upper and lower solutions again. From Remark \ref{remark3.5}, we know that there exists a lower solution of the above equation. Hence, it suffices to find an upper solution of \eqref{3.21.}. Since $\hat{\rho}\in PC(\theta)$, there is a positive solution $u\in C^\infty(M)$ of $-b_n\Delta_\theta u+\rho u=\hat{\rho} u^a$. If $\hat{\rho}_1\leq \hat{\rho}$, then $u$ is an upper solution of \eqref{3.21.}. Indeed, 
\begin{align}
-b_n\Delta_\theta u+\rho u-\hat{\rho}_1u^a=(-b_n\Delta_\theta u+\rho u-\hat{\rho}u^a)+(\hat{\rho}-\hat{\rho}_1)u^a\geq 0.
\end{align}
Hence we may get a positive solution of \eqref{3.21.}.
\qed
\begin{Remark}
If $\hat{\rho}\in PC(\theta)$ and $\hat{\rho}_1=\alpha\hat{\rho}$ for some constant $\alpha>0$, then $\hat{\rho}_1\in PC(\theta)$ regardless of the sign of $\lambda_1$. Indeed, since $\hat{\rho}\in PC(\theta)$, there is a positive solution $u\in C^\infty(M)$ of \eqref{2.10.}, then $\alpha^{-\frac{1}{a-1}}u$ is a solution of \eqref{3.21.}, so $\hat{\rho}_1\in PC(\theta)$.
\end{Remark}

Now we turn to the case $\lambda_1=0$. 

\begin{Proposition}
\label{theorem3.8.}
Let $(M^{2n+1},H,J,\theta)$ be a compact pseudo-Hermitian manifold. Then $0\in PC(\theta)$ if and only if $\lambda_1=0$.\end{Proposition}
\proof
Assume that $0\in PC(\theta)$, that is, there is a positive solution $u\in C^\infty(M)$ of $Lu=-b_n\Delta_\theta u+\rho u=0$. From
\begin{align}
\lambda_1\langle \psi, u\rangle_{L^2}=\langle L\psi, u\rangle_{L^2}=\langle \psi, Lu\rangle_{L^2}=0
\end{align}
 where $\psi$ is the positive eigenfunction of $\lambda_1$ of $L$, we deduce that $\lambda_1=0$. Conversely, if $\lambda_1=0$, then the associated eigenfunction $\psi$ realizes the zero Webster curvature, i.e., $0\in PC(\theta)$. 
\qed

Since by Proposition \ref{theorem3.8.} if $\lambda_1=0$ then one can always find a pseudo-Hermitian structure $\hat{\theta}\in [\theta]$ of zero Webster scalar curvature, we can without loss of generality restrict our attention to the case that $(M^{2n+1},H,J,\theta)$ already has a zero Webster scalar curvature $\rho\equiv0$.

\begin{Proposition}
Let $(M^{2n+1},H,J,\theta)$ be a compact pseudo-Hermitian manifold with Webster scalar curvature $\rho\equiv 0$. If  $0\not\equiv\hat{\rho}\in PC(\theta)$, then $\hat{\rho}$ must change sign on $M$ and $\int_M\hat{\rho}\ \Psi^\theta<0$.
\end{Proposition}
\proof
Since $\hat{\rho}\in PC(\theta)$, there is a positive solution $u\in C^\infty(M)$ such that
\begin{align}
-b_n\Delta_\theta u=\hat{\rho}u^{a},\label{3.27...}
\end{align}
where $b_n=2+\frac{2}{n}$, $a=1+\frac{2}{n}$. Hence,
\begin{align}
\int_M \hat{\rho}u^{a}\Psi^\theta=-b_n\int_M \Delta_\theta u\ \Psi^\theta=0.
\end{align}
Therefore, $\hat{\rho}$ must change sign on $M$ since $u>0$ and $\hat{\rho}\not\equiv 0$. Furthermore, multiplying \eqref{3.27...} by $u^{-a}$ and integrating by parts yield
\begin{align}
\int_M\hat{\rho}\Psi^\theta=-b_n\int_Mu^{-a}\Delta_\theta u\ \Psi^\theta=-ab_n\int_Mu^{-a-1}|\nabla^Hu|^2\Psi^\theta<0.
\end{align}
\qed

In the case $\lambda_1>0$, there is a substitute for Theorem \ref{theorem3.4.} as follows.

\begin{Proposition}
\label{theorem3.9.}
Let $(M^{2n+1},H,J,\theta)$ be a compact pseudo-Hermitian manifold. Then $\lambda_1>0$ if and only if there is a positive function $\hat{\rho}\in C^\infty(M)$ such that $\hat{\rho}\in PC(\theta)$.
\end{Proposition}
\proof
Let $\psi$ be the positive eigenfunction associated with $\lambda_1$ of $L$. If there is a positive function $\hat{\rho}\in C^\infty(M)$ such that $\hat{\rho}\in PC(\theta)$, then $Lu=\hat{\rho}u^a$ has a positive solution $u$, and 
\begin{align}
\lambda_1\langle \psi, u\rangle_{L^2}=\langle L\psi, u\rangle_{L^2}=\langle \psi, Lu\rangle_{L^2}=\langle \psi, \hat{\rho}u^a\rangle_{L^2}>0
\end{align}
which implies that $\lambda_1>0$. Conversely, if $\lambda_1>0$, then
\begin{align}
L\psi=\lambda_1\psi=(\lambda_1\psi^{1-a})\psi^a.
\end{align}
Pick $\hat{\rho}=\lambda_1\psi^{1-a}>0$, then $L\psi=\hat{\rho}\psi^a$, so $\hat{\rho}\in PC(\theta)$.
\qed
 
Combining Theorem \ref{theorem3.4.}, Proposition \ref{theorem3.8.}, Proposition \ref{theorem3.9.} with the definition of the CR Yamabe constant, we can give the proof of Theorem \ref{theorem3.10}.
\proof[\bf{Proof of Theorem \ref{theorem3.10}}]
First, we show the assertion $(A)$. Clearly, ``$(a1)\Leftrightarrow(a2)$" has already been proved in Theorem \ref{theorem3.4.}. The statement ``$(a3)\Leftrightarrow(a4)$" can be deduced from \eqref{2.15} and Theorem 3.4 of \cite{[JL1]}. It remains to prove the statement ``$(a2)\Leftrightarrow(a3)$". The case ``$(a2)\Rightarrow(a3)$" is trivial. Let us consider ``$(a3)\Rightarrow(a2)$". Assume that ``$\hat{\rho}_0\in \{\hat{\rho}\in C^\infty(M):\hat{\rho}<0\} \cap PC(\theta)$". It follows from Theorem \ref{theorem3.4.} that $\lambda_1<0$. In terms of the equivalence ``$(a1)\Leftrightarrow(a2)$", we see that $(a2)$ holds.

Next, we consider the assertion $(B)$. Obviously, ``$(b1)\Leftrightarrow(b2)$" has been proved in Proposition \ref{theorem3.8.}. ``$(b3)\Rightarrow(b2)$" follows from Theorem 3.4 of \cite{[JL1]}. Now let us consider the case ``$(b2)\Rightarrow(b3)$". By \eqref{2.15}, we have $Y_M(\theta)\leq 0$. If $Y_M(\theta)<0$, then $\lambda_1<0$ by the results in $(A)$ of this theorem. However, according to ``$(b1)\Leftrightarrow(b2) $", $0\in PC(\theta)$ implies $\lambda_1=0$, which leads to a contradiction. Therefore, $Y_M(\theta)=0$.

Finally, we treat the assertion $(C)$. It is obvious that ``$(c1)\Leftrightarrow(c2)$" can be obtained by Proposition \ref{theorem3.9.}. ``$(c3)\Rightarrow(c2)$" may be deduced by the solvability of CR Yamabe problem (cf. \cite{[Ga1]}, \cite{[GY]}, \cite{[JL1]}). For the case ``$(c2)\Rightarrow(c3)$", using the equivalence ``$(c1)\Leftrightarrow(c2)$", we have $\lambda_1>0$. Combining ``$(a1)\Leftrightarrow(a4)$" and ``$(b1)\Leftrightarrow(b3)$" yields $Y_M(\theta)>0$.
\qed

From Theorem \ref{theorem3.10}, we can easily see that
\begin{Corollary}\label{corollary3.11}
Let $(M^{2n+1},H,J,\theta)$ be a compact pseudo-Hermitian manifold. Then $\lambda_1$ and $Y_M(\theta)$ have the same sign or are both zero, and thus the sign of $\lambda_1$ is a CR invariant. 
\end{Corollary}
\begin{Remark}
In \cite{[Wa]}, it is easy to see that Corollary \ref{corollary3.11} is implicit in his argument, although it is not clearly pointed out. Indeed, it was proved directly by using the transformation law $(3.1)$ in \cite{[JL1]} under the CR pointwise conformal deformations and the solvability of CR Yamabe problem.
\end{Remark}

\section{CR conformally equivalent deformations with prescribed Webster scalar curvature}
In this section, we will determine the set $CE(\theta)$ on a compact pseudo-Hermitian manifold $(M^{2n+1},H,J,\theta)$ with $\lambda_1<0$, $\lambda_1=0$ and $\lambda_1>0$ respectively.

Let us consider a second-order quasilinear degenerate elliptic differential operator $T$ on the compact pseudo-Hermitian manifold $(M^{2n+1}, H, J, \theta)$:
\begin{align}
 Tu=u^{-a}\left(-b_n\Delta_\theta u+\rho u \right),\label{3.1}
\end{align}
where $a=1+\frac{2}{n}$, $b_n=2+\frac{2}{n}$ and $\rho\in C^\infty(M)$. The linearization of $T$ at a given positive function $u_0\in C^\infty(M)$ is a second-order linear degenerate elliptic differential operator
 \begin{align}
T'(u_0)v&=\left.\frac{d}{dt}\right |_{t=0}T(u_0+tv)\notag\\
          &=b_nu_0^{-a}\left\{-\Delta_\theta v+\left(a\frac{\Delta_\theta u_0}{u_0}+\frac{1-a}{b_n}\rho\right) v\right\},\label{4.2.}
 \end{align}
where $v\in S^p_2(M)$. Set
\begin{align}
A(u_0)v=-\Delta_\theta v+\left( a\frac{\Delta_\theta u_0}{u_0}+\frac{1-a}{b_n}\rho\right)v,\label{3.14}
\end{align}
which is a linear self-adjoint degenerate elliptic operator with $\ker T'(u_0)=\ker A(u_0)$. 
 
\begin{Lemma}\label{lemma3.1}
Let $(M^{2n+1},H,J,\theta)$ be a compact pseudo-Hermitian manifold. Let $L_3: S^p_2(M)\rightarrow L^p(M)$ be the operator defined as in \eqref{4.2.} with $0<u_0\in C^\infty(M)$ and $\rho\in C^\infty(M)$. Assume that $p>2n+1$. If $\ker L_3=0$, then
\begin{align}
\|v\|_{S^p_2(M)}\leq C\|L_3v\|_{L^p(M)}\label{3.5}
\end{align}
for any $v\in S^p_2(M)$, where $C$ is a positive constant independent of $v$. Therefore, the operator $L_3: S^p_2(M)\rightarrow L^p(M)$ is bijective with a continuous inverse.
\end{Lemma}
\proof
Since $M$ is compact, there exists a constant $C_1>0$ such that 
\begin{align}
\|v\|_{S^p_2(M)}\leq C_1\left(\|L_3v\|_{L^p(M)}+\|v\|_{L^p(M)} \right)\label{3.6}
\end{align}
for any $v\in S^p_2(M)$, which can be deduced from the $L^p$ interior regularity results for $\Delta_\theta$ (cf. Theorem 18 of \cite{[RS]}) by using a partition of unity on $M$. In order to get \eqref{3.5}, it is sufficient to prove that 
\begin{align}
\|v\|_{L^p(M)}\leq C_2\left\|L_3v\right\|_{L^p(M)}
\end{align}
 for any $v\in S^p_2(M)$, where $C_2>0$ is a constant independent of $v$. If not, there is a sequence $\{v_n\}\subset S^p_2(M)$ such that $\|v_n\|_{L^p(M)}=1$ but $\left\|L_3v_n\right\|_{L^p(M)}\rightarrow 0$ as $n\rightarrow +\infty$. Then by \eqref{3.6}, we have
 \begin{align}
 \|v_n\|_{S^p_2(M)}\leq C_3,
 \end{align}
 where the constant $C_3>0$ is independent of $n$. Using the compactly embedding theorem $S^p_2(M)\subset W^{1,p}(M)\subset\subset C^0(M)$ (cf. Theorem 19.1 of \cite{[FS]}) yields that there exists a subsequence $v_{n_k}$ of $v_n$ and a function $v\in C^0(M)$ such that $\lim_{k\rightarrow +\infty}\|v_{n_k}-v\|_{C^0(M)}=0$, and thus $\lim_{k\rightarrow +\infty}\|v_{n_k}-v\|_{L^p(M)}=0$, where $W^{1,p}(M)$ is the classical Sobolev space with $p>2n+1$. According to \eqref{3.6} and the triangle inequality, we have
 \begin{align}
 \|v_{n_i}-v_{n_j}\|_{S^p_2(M)}\leq C_1&\left( \|L_3v_{n_i}\|_{L^p(M)}+\|L_3v_{n_j}\|_{L^p(M)}\right.\notag\\
 &\left.+\|v_{n_i}-v\|_{L^p(M)}+\|v-v_{n_j}\|_{L^p(M)}\right)\rightarrow 0
 \end{align}
 as $i,j\rightarrow \infty$, i.e., $\{v_{n_k}\}$ is a Cauchy sequence in $S^p_2(M)$, so $\lim_{k\rightarrow +\infty} \|v_{n_k}-v\|_{S^p_2(M)}=0$. By the  continuity of the operator $L_3: S^p_2(M)\rightarrow L^p(M)$, we obtain that 
 \begin{align}
 \|L_3v\|_{L^p(M)}=\lim_{k\rightarrow +\infty}\|L_3v_{n_k}\|_{L^p(M)}=0.
 \end{align}
Hence, $L_3v=0$ in $L^p(M)$. By $\ker L_3=0$, we get $v=0$ in $S^p_2(M)$. However, from $\|v_{n_k}\|_{L^p(M)}=1$ for any $k$, we deduce that $\|v\|_{L^p(M)}=1$, which leads to a contradiction. 

Now we go to prove the last conclusion of this lemma, namely, $L_3: S^p_2(M)\rightarrow L^p(M)$ is bijective with a continuous inverse. In fact, the condition $\ker L_3=0$ implies the injectivity. So it is sufficient to show the existence of $L_3v=f$ for any $f\in L^p(M)$. According to the fact that $C^\infty(M)$ is dense in $L^p(M)$, there is a sequence $\{f_j\}\subset C^\infty(M)$ such that $\lim_{j\rightarrow \infty}\|f_j-f\|_{L^p(M)}=0$. In terms of  Lemma \ref{lemma4.1} and regularity results in \cite{[Xu-2]}, there exists $v_j\in C^\infty(M)$ such that $L_3v_j=f_j$. Using \eqref{3.5} and $\lim_{j\rightarrow \infty}\|f_j-f\|_{L^p(M)}=0$ yields that 
\begin{align}
\|v_i-v_j\|_{S^p_2(M)}\leq C\|L_3v_i-L_3v_j\|_{L^p(M)}=C\|f_i-f_j\|_{L^p(M)}\rightarrow 0
\end{align}
 as $i,j\rightarrow +\infty$, so $\{v_i\}$ is a Cauchy sequence in $S^p_2(M)$. Hence, $v_j\rightarrow v$ in $S^p_2(M)$ due to the completeness of $(S^p_2(M), \|\cdot\|_{S^p_2(M)})$. Since $L_3: S^p_2(M)\rightarrow L^p(M)$ is a continuous map, 
 \begin{align}
 \|L_3v-f\|_{L^p(M)}&\leq  \|L_3v-L_3v_j\|_{L^p(M)}+ \|L_3v_j-f\|_{L^p(M)}\notag\\
 &=\|L_3v-L_3v_j\|_{L^p(M)}+ \|f_j-f\|_{L^p(M)}\rightarrow 0
 \end{align}
 as $j\rightarrow +\infty$. Consequently, $L_3v=f$. Hence, $L_3: S^p_2(M)\rightarrow L^p(M)$ is a bijective continuous linear map, and thus $L_3$ has a continuous inverse map due to the open mapping theorem of Banach.
\qed

Using the inverse function theorem for Banach spaces and the regularity results for degenerate elliptic equations (cf. \cite{[Xu-2]}), we have the following theorem.
\begin{Theorem}\label{theorem3.2}
Let $(M^{2n+1},H,J,\theta)$ be a compact pseudo-Hermitian manifold, and let $T$ be the operator defined as in \eqref{3.1} with $\rho\in C^\infty(M)$. Assume that $0<u_0\in C^\infty(M),\ p>2n+1$. If the linearization $T'(u_0):S^p_2(M)\rightarrow L^p(M)$ is injective (and thus is invertible), then there exists a constant $\eta>0$ such that for any function $f\in C^\infty(M)$ with $\|f-T(u_0)\|_{L^p(M)}<\eta$, there is a positive function $u\in C^\infty(M)$ satisfying $T(u)=f$.
\end{Theorem}

According to Theorem 4.10 of Chapter five of \cite{[Ka]}, the spectrum of the self-adjoint operator $A(u)$ depends continuously on $u$. Furthermore, from the proof of Lemma \ref{lemma4.1} (2), we know that the resolvent of the self-adjoint operator $A(u)$ is compact for any positive function $u\in C^\infty(M)$. Therefore, the spectrum of $A(u)$ is discrete, and every eigenvalue is of finite multiplicity. In addition, given a function $z\in C^\infty(M)$, the self-adjoint operator $A(u+tz)$ depends analytically on $t$ for $|t|$ small enough, hence so do the eigenvalues and eigenfunctions of $A(u+tz)$ (cf. \cite{[KMR]}). After a process similar to Theorem 4.5 and Lemma 4.6 in \cite{[KW2]}, we have the following perturbation theorem for $T'$.

\begin{Theorem}\label{theorem3.4}
 The second-order linear degenerate elliptic operator $T'(u): S^p_2(M)\rightarrow L^p(M)$ is bijective on an open dense subset of  the set $\{u\in C^\infty(M):\ u>0\}$.
\end{Theorem}

For our purpose in this section, we need the following approximation theorem (cf. Theorem 2.1 of \cite{[KW2]}).
\begin{Theorem}\label{theorem4.4.}
Let $N$ be a connected manifold with dimension $n\geq 2$ and let $f\in C(N)\cap L^p(N)$. Then a function $g\in L^p(N)$ is in the $L^p-$closure of $O_f$ if and only if $\inf_N f\leq g(x)\leq \sup_N f$ for almost all $x\in N$. Here $O_f$ is the orbit of $f$ under the group of diffeomorphism of $N$.
\end{Theorem}

Making use of the above three theorems yields the key lemma as follows.

\begin{Lemma}\label{lemma3.5}
For a given smooth function $\hat{\rho}$ on a compact pseudo-Hermitian manifold $(M^{2n+1},H,J,\theta)$ with the Webster scalar curvature $\rho$, if $\min_{M} \hat{\rho}<C\rho<\max_M\hat{\rho}$ for some constant $C>0$, then $\hat{\rho}\in CE(\theta)$.
\end{Lemma}
\proof
Let $u_0\equiv1$. According to Theorem \ref{theorem3.4}, for any $\epsilon >0$, there exists a smooth function $u_1$ so close to $u_0$ that 
\begin{align}
\|T(u_1)-\rho\|_{\infty}=\|T(u_1)-T(u_0)\|_{\infty}<\epsilon,
\end{align}
and $T'(u_1)$ is invertible. Picking $\epsilon$ sufficiently small and using the assumption $C^{-1}\min_{M} \hat{\rho}<\rho<C^{-1}\max_M\hat{\rho}$ yield
\begin{align}
C^{-1}\min_{M} \hat{\rho}<T(u_1)<C^{-1}\max_M\hat{\rho}.
\end{align}
By Theorem \ref{theorem4.4.}, we obtain that for any $\eta>0$, there exists a diffeomorphism $\Phi$ of $M$ such that
\begin{align}
\|C^{-1}\hat{\rho}\circ \Phi-T(u_1)\|_{L^p(M)}<\eta.
\end{align}
Making use of Theorem \ref{theorem3.2}, we get that there exists a positive solution $u\in C^\infty(M)$ of $T(u)=C^{-1}\hat{\rho}\circ \Phi$. Set $v=C^{-\frac{1}{a-1}}u$, then $v$ is a positive solution of $T(v)=\hat{\rho}\circ\Phi$. Therefore, $\hat{\rho}\in CE(\theta)$.
\qed

In terms of the above key lemma, we can give the proof of Theorem \ref{theorem1.2}. 
\proof[\bf{Proof of Theorem \ref{theorem1.2}}]
$(1)$ Since $\lambda_1<0$ and Theorem \ref{theorem3.4.}, there is a pseudo-Hermitian structure $\theta_1\in [\theta]$ such that the corresponding Webster scalar curvature equals -1. If $\hat{\rho}$ is negative somewhere, then $\min_M \hat{\rho}<-C<\max_M \hat{\rho}$ for some constant $C>0$. According to Lemma \ref{lemma3.5}, we have $\hat{\rho}\in CE(\theta_1)=CE(\theta)$. Conversely, if $\hat{\rho}\in CE(\theta)$, then there exists a diffeomorphism $\Phi$ of $M$ and a positive function $u\in C^\infty(M)$ such that $Lu=(\hat{\rho}\circ \Phi) u^a$. Let $\psi$ be the positive eigenfunction associated with the eigenvalue $\lambda_1$ of $L$, then
\begin{align}
0>\lambda_1\langle \psi, u\rangle_{L^2}=\langle L\psi, u\rangle_{L^2}=\langle \psi, Lu\rangle_{L^2}=\langle \psi, (\hat{\rho}\circ \Phi) u^a\rangle_{L^2}.\label{3.35}
\end{align}
Consequently, $\hat{\rho}$ must be negative somewhere on $M$.

$(2)$ From Proposition \ref{theorem3.8.}, it follows that $0\in PC(\theta)\subset CE(\theta)$. Similar to part (1) of this theorem, we can get the conclusion easily.

$(3)$ According to Theorem \ref{theorem3.10} and the results about CR Yamabe problem (cf. \cite{[Ga1]}, \cite{[GY]}, \cite{[JL1]}), $\lambda_1>0$ implies that there is a positive constant $\rho_1\in PC(\theta)$. By an argument similar to part (1) of this theorem, we can obtain that $\hat{\rho}\in CE(\theta)$ if and only if $\hat{\rho}$ is positive somewhere. 
\qed

Before the end of this section, we point out that the sign of $\lambda_1$ is invariant under CR conformally equivalent deformations.
\begin{Theorem}
Let $(M^{2n+1},H,J,\theta)$ be a compact pseudo-Hermitian manifold. If the structures $(\hat{H},\hat{J},\hat{\theta})$ and $(H,J,\theta)$ are CR conformally equivalent, then $\lambda_1(\theta)$ and $\lambda_1(\hat{\theta})$ have the same sign or are both zero.
\end{Theorem}
\proof
Since $(\hat{H},\hat{J},\hat{\theta})$ and $(H,J,\theta)$ are CR conformally equivalent, then there is a map $\Phi\in \text{Diff}(M)$ and $0<u\in C^\infty(M)$ such that
\begin{align}
\Phi^*\hat{\theta}=u^{\frac{2}{n}}\theta,\ \ \ \hat{H}=d\Phi(H),\ \ \ \hat{J}=d\Phi\circ J\circ (d\Phi)^{-1}.
\end{align}
Set $\tilde{\theta}=u^{\frac{2}{n}}\theta$, hence $\Phi: (M,H,J,\tilde{\theta})\rightarrow (M, \hat{H}, \hat{J}, \hat{\theta})$ is a CR isomorphism with $\Phi^*\hat{\theta}=\tilde{\theta}$. Therefore, $\lambda_1(\hat{\theta})=\lambda_1(\tilde{\theta})$. From Corollary \ref{corollary3.11}, it follows that $\lambda_1(\hat{\theta})$ and $\lambda_1(\theta)$ have the same sign or are both zero.
\qed

\section*{Appendix}
In this section, we give an alternative proof of Theorem \ref{theorem3.4.} by following the spirit of \cite{[KW2]}. Although the expressions are different, the following theorem and Theorem \ref{theorem3.4.} are completely equivalent. 

\begin{Theorem}
Let $(M^{2n+1},H,J,\theta)$ be a compact pseudo-Hermitian manifold. Then $S=\{f\in C^\infty(M):\ f<0\}\subset \text{Im}\ T$ if and only if $\lambda_1<0$, where $T$ is as in \eqref{3.1} with $\rho\in C^\infty(M)$, and $\text{Im}\ T=\{Tu: 0<u\in S^p_2(M)\}$.
\end{Theorem}
\proof
If $S\subset  \text{Im}\ T$, then $-1\in \text{Im}\ T$, i.e., $T(u)=-1$ for some $C^\infty(M)$ function $u>0$. Thus, $Lu=-u^a$. Let $\psi$ be the positive eigenfunction of $L$ with respect to $\lambda_1$, namely, $L\psi=\lambda_1\psi$. So we have the following
\begin{align}
\lambda_1\langle \psi, u\rangle_{L^2}=\langle L \psi, u\rangle_{L^2}=\langle \psi, Lu\rangle_{L^2}=-\langle \psi, u^a\rangle_{L^2}<0
\end{align}
 which implies $\lambda_1<0$. 
 
 For the converse, assume that $\lambda_1<0$. Set $K=S\cap \text{Im}\ T$. From $L\psi=\lambda_1\psi$, we have $T(\psi)=\lambda_1\psi^{1-a}<0$, thus $\lambda_1\psi^{1-a}\in K$. Consequently, $K$ is nonempty. Clearly, $S$ is connected. In order to prove $K=S$, it is sufficient to prove $K$ is a both open and closed subset in $S$. For openness, we will show that for any $u>0$, $T(u)\in K$ implies $\ker T'(u)=\ker A(u)=0$, where $A(u)$ is defined by \eqref{3.14}, and thus $K$ is open subset of $S$ in terms of Theorem \ref{theorem3.2}. Let $\mu_1$ be the first eigenvalue of $A(u)$ and $\phi$ be the corresponding positive eigenfunction. By \eqref{3.14}, we have 
 \begin{align}
 \mu_1\langle \phi, u\rangle_{L^2}=\langle A(u)\phi, u\rangle_{L^2}=\langle \phi, A(u)u\rangle_{L^2}=\frac{1-a}{b_n}\langle \phi, T(u)u^a \rangle_{L^2}>0
 \end{align}
 which gives $\mu_1>0$, and so $\ker A(u)=0$. For closeness, we assume that $f_j\in K$ and $f_j\xrightarrow{C^0} f\in S$, we need prove $f\in K$, i.e., there is $u\in S^p_2(M)$ such that $T(u)=f$. Since $f_j\in K$, there exists a function $0<u_j\in C^\infty(M)$ satisfying $T(u_j)=f_j$. Let $w_j=\log{\frac{u_j}{\psi}}$ where $\psi$ is the positive eigenfunction associated with the first eigenvalue $\lambda_1$ of $L$. Then $w_j$ satisfies 
 \begin{align}
 -b_n\Delta_\theta w_j-b_n \nabla^H w_j\cdot\left(\nabla^H w_j+2\frac{\nabla^H\psi}{\psi} \right)=-\lambda_1+f_j\psi^{a-1}e^{(a-1)w_j}
 \end{align}
 where $\cdot$ is the inner product induced by the Webster metric $g_\theta$. Considering the maximum and minimum of $w_j$ and using the classical maximum principle, it is easy to show that there are two constants $m_1, m_2>0$ independent of $j$ such that $0<m_1\leq u_j\leq m_2$.  Hence, applying Lemma \ref{lemma3.1} to the operator $L_4=-\Delta_\theta+ id$, we have
\begin{align}
\|u_j\|_{S^p_2(M)}\leq C\|L_4u_j\|_{L^p(M)}=C\left\|\frac{1}{b_n}(f_ju_j^a-\rho u_j)+u_j\right\|_{L^p(M)}\leq \hat{C}
\end{align}
 where $C, \hat{C}$ are constants independent of $j$. Using the compactly embedding theorem $S^p_2(M)\subset W^{1,p}(M)\subset\subset C^0(M)$ ( $p>2n+1$, cf. Theorem 19.1 of \cite{[FS]}), there exists a subsequence $\{u_{j_k}\}$ such that $u_{j_k}\xrightarrow{C^0} u$ as $k\rightarrow +\infty$, where $u>0$ in $M$ since $0<m_1\leq u_{j_k}\leq m_2$. Moreover, the subsequence $\{u_{j_k}\}$ is a Cauchy sequence in $S^p_2(M)$, because
 \begin{align}
 \|u_{j_k}-u_{j_l}\|_{S^p_2(M)}&\leq C \|L_4(u_{j_k}-u_{j_l})\|_{L^p(M)}\notag\\
       						&=C\left\|(u_{j_k}-u_{j_l})+\frac{1}{b_n}(f_{j_k}u_{j_k}^a-f_{j_l}u_{j_l}^a+\rho u_{j_l}-\rho u_{j_k})\right\|_{L^p(M)}\notag\\
						&\rightarrow 0,
 \end{align}
 as $k,l \rightarrow +\infty$. Therefore, $u_{j_k}\rightarrow u$ in $S^p_2(M)$ as $k\rightarrow +\infty$. Let $k\rightarrow\infty$ in $T(u_{j_k})=f_{j_k}$, by the continuity of $T: S^p_2(M)\rightarrow L^p(M)$, we obtain $T(u)=f$, so $f\in K$. 

\qed

\bigskip

Yuxin Dong \ \ 

School of Mathematical Sciences

Fudan University

Shanghai, 200433, P. R. China

yxdong@fudan.edu.cn \ \ \ 

\bigskip

Yibin Ren

College of Mathematics and Computer Science

Zhejiang Normal University

Jinhua, 321004, Zhejiang, P.R. China

allenryb@outlook.com

\bigskip

Weike Yu

School of Mathematical Sciences

Fudan University

Shanghai, 200433, P. R. China

wkyu2018@outlook.com

\bigskip

\end{document}